\documentclass[a4paper,11pt]{article}

\usepackage{graphicx}   % Graphics
\usepackage{euscript}  % Euler fonts
\usepackage{amsfonts}  % Euler fonts
\usepackage{amssymb}
\usepackage{latexsym}
\usepackage{amscd}
\usepackage{harvard}
\newcommand{\field}[1]{\mathbb{#1}} 
\newtheorem{thm}{Theorem}[section] 
\newtheorem{mainthmm}{Main Theorem}
\newtheorem{oldthm}{Theorem}

\newtheorem{lem}[thm]{Lemma} 
\newtheorem{prop}[thm]{Proposition}

\newtheorem{defn}[thm]{Definition}
 
\newtheorem{rmk}[thm]{Remark}

\def\C{\mathbb{C}}

\def\be{\beta}

\def\ups{\upsilon}
\def\ga{\gamma}

\def\proof{{\it Proof. \ }}
\def\QED{{\flushright $\Box$}}

\def\gmu{\mu}

\def\fR{\field{R}}

\def\cA{\mathcal{A}}
\def\cB{\mathcal{B}}

\def\cD{\mathcal{D}}

\def\cI{\mathcal{I}}

\def\cN{\mathcal{N}}
\def\cM{\mathcal{M}}

\def\cR{\mathcal{R}}

\def\cS{\mathcal{S}}

\def\cZ{\mathcal{Z}}

\begin{document}

\title{Period Doubling Renormalization for Area-Preserving Maps and Mild Computer Assistance in Contraction Mapping Principle}

\author{Denis G. Gaidashev \\
\\
Department of Mathematics, \\
Uppsala University, Box 480, \\
751 06 Uppsala, Sweden, \\
{\tt gaidash@math.uu.se}}

\maketitle

\begin{abstract}

A universal period doubling cascade analogous to the famous Feigenbaum-Coullet-Tresser period doubling has been observed in area-preserving maps of ${\fR}^2$. Existence of the ``universal'' map with orbits of all binary periods has been proved via a renormalization approach in \cite{EKW2} and \cite{GJM}. These proofs  use ``hard'' computer assistance.

%A  renormalization approach has been used in a ``hard'' computer-assisted proof of existence of an area-preserving map with orbits of all binary periods in \cite{EKW2}. As it is the case with all non-trivial universality problems in non-dissipative systems in dimensions more than one, no analytic proof of this period doubling universality exists to date. 

In this paper we attempt to reduce computer assistance in the argument, and present a mild computer aided proof of the analyticity and compactness of the renormalization operator  in a neighborhood of a renormalization fixed point: that is a proof that does not use generalizations of interval arithmetics to functional spaces ---  but rather relies on interval arithmetics on real numbers only to estimate otherwise explicit expressions. The proof relies on several instance of the Contraction Mapping Principle,  which is, again, verified via mild computer assistance.

\end{abstract}

\newpage

\setcounter{page}{1}

\section{Introduction: Period doubling universality for area-preserving maps}

Following the pioneering discovery of the Feigenbaum-Coullet-Tresser period doubling universality in unimodal maps \cite{Fei1,Fei2,TC} universality has been demonstrated to be a rather generic phenomenon in dynamics.

To prove universality one usually introduces a {\it renormalization} operator in a functional space, and demonstrates that this operator has a hyperbolic fixed point.

Such renormalization approach to universality has been very successful in one-dimensional dynamics, and has led to explanation of universality in unimodal maps \cite{Eps1,Eps2,Sul,Lyu}, critical circle maps \cite{dF1,dF2,Ya1,Ya2} and holomorphic maps with a Siegel disks \cite{McM,Ya3,GY}. Universality has been also abundantly observed in higher dimensions, in particular, in two and more dimensional dissipative systems (cf.  \cite{CEK1}, \cite{Spa}), in area-preserving maps, both as the period doubling universality \cite{DP,Hel,BCGG,CEK2,EKW1,EKW2,GK2,GJ1,GJ2,GJM}, and as the universality associated with the break-up of invariant surfaces  \cite{SK,McK1,McK2,ME},  and in Hamiltonian flows \cite{ED,AK,AKW,Koch1,Koch2,Koch3,GK1,Gai1,Kocic,KLDM}.

In this paper we will consider a period doubling universality for area-preserving maps of the plane ---  an analogue of Feigenbaum-Coullet-Tresser universality in higher dimensions. 

An infinite period doubling cascade in families of area-preserving maps was observed by several authors in  early 80's \cite{DP,Hel,BCGG,Bou,CEK2}. A typical period doubling scenario can be illustrated with the area-preserving H\' enon family (cf. \cite{Bou}) :
$$ H_a(x,y)=(-y +1 - a x^2, x).$$

Maps in this family posses a fixed point $((-1+\sqrt{1+a})/a,(-1+\sqrt{1+a})/a) $ which is stable for $-1 < a < 3$. When $a_1=3$ this fixed point becomes unstable, at the same time an orbit of period two is born with $H_a(x_\pm,x_\mp)=(x_\mp,x_\pm)$, $x_\pm= (1\pm \sqrt{a-3})/a$. This orbit, in turn, becomes unstable at $a_2=4$, giving birth to a period $4$ stable orbit. Generally, there  exists a sequence of parameter values $a_k$, at which the orbit of period $2^{k-1}$ turns unstable, while at the same time a stable orbit of period $2^k$ is born. The parameter values $a_k$ accumulate on some $a_\infty$. The crucial observation is that the accumulation rate
$$\lim_{k \rightarrow \infty}{a_k-a_{k-1} \over  a_{k+1}-a_k } = 8.721...$$
is universal for a large class of families, not necessarily H\'enon.

Furthermore, the $2^k$ periodic orbits scale asymptotically with two scaling parameters
$$
\lambda=-0.249...,\quad \gmu=0.061...
$$

To explain how orbits scale with $\lambda$ and $\mu$ we will follow \cite{Bou}. Consider an interval $(a_k,a_{k+1})$ of parameter values in a `''typical'' family $F_a$. For any value $\alpha \in (a_k,a_{k+1})$ the map $F_\alpha$ posses a stable periodic orbit of period $2^k$. We fix some $\alpha_k$ within the interval $(a_k,a_{k+1})$ in some consistent way; for instance, by requiring that the restriction of $F^{2^k}_{\alpha_k}$ to a neighborhood of a stable periodic point in the $2^{k}$-periodic orbit is conjugate via a diffeomorphism $H_k$ to a rotation with some fixed rotation number $r$.  Let $p'_k$ be some unstable periodic point in the $2^{k-1}$-periodic orbit, and let $p_k$ be the further (from $p'_k$) of the two stable $2^{k}$-periodic points that bifurcated from $p'_k$.  Then,
$${1 \over \lambda}=-\lim_{k \rightarrow \infty}{ |p_k-p'_k | \over |p_{k+1}-p'_{k+1} |},\quad {\lambda \over \mu}=-\lim_{k \rightarrow \infty}{ \rho_k \over \rho_{k+1}},$$
%\lim_{k->\infty}{ c_k \over c_{k+1}}&=& {1 \over \lambda^2}.$$
where $\rho_k$ is the ratio of the eigenvalues of $D H_k(p_k)$.

This universality can be explained rigorously if one shows that the {\it renormalization} operator
\begin{equation}\label{Ren}
R[F]=\Lambda^{-1}_F \circ F \circ F \circ \Lambda_F,
\end{equation}
where $\Lambda_F$ is some $F$-dependent coordinate transformation, has a fixed point, and the derivative of this operator is hyperbolic at this fixed point.

It has been argued in \cite{CEK2}  that $\Lambda_F$ is a diagonal linear transformation. Furthermore, such $\Lambda_F$ has been used in \cite{EKW1} and \cite{EKW2} in a computer assisted proof of existence of the renormalization fixed point. This is the strongest result to date concerning the existence of the renormalization  fixed point.

The problem has been further studied in \cite{GK2}, where the authors have demonstrated that  the period doubling renormalization fixed point of \cite{EKW2} is ``almost'' one dimensional, in the sense that it is close to the following H\'enon-like map:
$$H_*(x,u)=(\phi(x)-u,x-\phi(\phi(x)-u )),$$
where $\phi$ solves
$$\phi(x)={2 \over \lambda} \phi(\phi(\lambda x))-x.$$ 

\cite{GK2} also suggests an analytic approach to the full period doubling problem for area-preserving maps based on its proximity to the one dimensional.

Infinitely renormalizable maps, that  is maps on the renormalization stable manifold, have been studied in \cite{GJ2} and \cite{GJM}, where it has been shown that such maps (in a codimension $1$ stable submanifold) posses  stable Cantor sets (sets on which the Lyapunov exponents are equal to zero) on which dynamics is rigid:  the stable dynamics for two different infinitely renormalizable maps  is conjugated by a $C^{1+\alpha}$ transformation (that is, a  differentiable transformation whose derivative is H\"older with exponent $\alpha$). This is in a stark contrast to the situation in dissipative H\'enon-like maps, where universality of Cantor sets (around the so called ``tip'') coexists with non-rigidity of the conjugacy of the dynamics on these sets (see \cite{dCLM,LM,HLM}).

The computer-assisted proof of \cite{EKW2} can be described as an application of the Contraction Mapping Principle for a Newton map for the operator $(\ref{Ren})$ in some neighborhood in an appropriately chosen Banach space. The proof is similar in spirit to the proofs \cite{Lan1,Lan2} of the Feigenbaum-Coullet-Tresser universality by Oscar Lanford, and uses bounds for analytic functions (not necessarily polynomials) via the so called ``standard sets'': loosely, generalizations of intervals in $\field{R}$ to infinite-dimensional functional spaces. Such computer assisted proofs came to be called ``hard'' within the computational community.  The approach turned out  to be quite powerful, albeit rather specialized and time and effort consuming.

In this paper we study  analyticity and compactness of the renormalization operator $(\ref{Ren})$  using, almost exclusively, the interval arithmetics in $\field{R}$. We obtain a bound on the neighborhood of  analyticity of the operator.  Specifically, we identify  a low order polynomial approximation for the ``generating function'' of the fixed point map  of $(\ref{Ren})$, and estimate the radius of the neighborhood of analyticity of the renormalization operator, centered on this polynomial approximation in an appropriate Banach space of functions analytic on a certain domain in $\field{R}^2$.

Furthermore, we demonstrate that renormalizations of functions from this neighborhood of analyticity, are themselves functions, analytic on a larger domain in $\field{R}^2$, which implies compactness of the renormalization operator.

We achieve these two goals using only what can be termed as ``mild'' computer assistance: computer estimates on arithmetic operations on reals and on the square root function on reals, as well as norm estimates of low order (quadratic) polynomials.  It is also our hope that the reader can find of interest several little ``tricks'' used in the paper to get bounds on solutions of functional equations.

\setcounter{page}{1}

\section{A renormalization operator on generating functions}

An ``area-preserving map'' will mean an exact symplectic diffeomorphism of a subset of ${\fR}^2$ onto its image.

Recall, that an area-preserving map $F: (x,u) \mapsto (y,v)$ can be uniquely specified by its generating function $\cS$:
\begin{equation}\label{gen_func}
\left( x \atop -\cS_1(x,y) \right) {{ \mbox{{\small \it  F}} \atop \mapsto} \atop \phantom{\mbox{\tiny .}}} \left( y \atop \cS_2(x,y) \right), \quad \cS_i \equiv \partial_i \cS,
\end{equation}
if the equation $u=-\cS_1(x,y)$ has the unique solution $y(x,u)$.

Furthermore, we will assume that $F$ is reversible, that is 
$$T \circ F \circ T=F^{-1}, \quad {\rm where} \quad T(x,u)=(x,-u).$$

For such maps it follows from $(\ref{gen_func})$ that 
$$\cS_1(y,x)=\cS_2(x,y) \equiv s(x,y),$$
and
$$
\nonumber \left({x  \atop  -s(y,x)} \right) {\mbox{\small \it F} \atop \mapsto} \left({y \atop s(x,y) }\right).
$$

It is this ``little'' $s$ that will be referred to below as ``the generating function''.

Applying a reversible $F$ twice we get
$$
 \left({x'  \atop  -s(z',x')} \right) {\mbox{\small \it F} \atop \mapsto} \left({z' \atop s(x',z')} \right)=
\left({z'  \atop  -s(y',z')} \right) {\mbox{\small \it F} \atop \mapsto} \left({y'\atop  s(z',y')} \right).
$$

It has been argued in \cite{CEK2}  that 

$$\Lambda_F(x,u)=(\lambda x, \gmu u),$$
where $\lambda$ and $\mu$ are some real parameters, dependent on $F$.

We therefore set  $(x',y')=(\lambda x,  \lambda y)$, $z'(\lambda x, \lambda y)= z(x,y)$ to obtain:

\begin{equation}\label{doubling}
\left(\!{x  \atop  -{ 1 \over \gmu } s(z,\lambda x)} \!\right) \!{\mbox{  {\small \it $\Lambda$}  }\atop \mapsto} \!\left(\!{\lambda x  \atop  -s(z,\lambda x)} \!\right) \!{\mbox{\small \it $F \circ F$} \atop \mapsto }\!\left(\!{\lambda y \atop s(z,\lambda y)}\! \right)   {\mbox{\small\it $\Lambda^{-1}$} \atop \mapsto} \left(\!{y \atop {1 \over \gmu } s(z,\lambda y) }\!\right),
\end{equation}
where $z(x,y)$ solves
\begin{equation}\label{midpoint}
s(\lambda x, z(x,y))+s(\lambda y, z(x,y))=0.
\end{equation}

If the solution of $(\ref{midpoint})$ is unique, then $z(x,y)=z(y,x)$, and it follows from $(\ref{doubling})$ that the generating function of the renormalized $F$ is given by 
\begin{equation}
\tilde{s}(x,y)=\gmu^{-1} s(z(x,y),\lambda y).
\end{equation}

Furthermore, it is possible to fix some normalization conditions for $\tilde{s}$ and $z$ which serve to determine scalings $\lambda$ and $\gmu$ as functions of $s$. Notice, that the normalization
$$s(1,0)=0$$ 
is reproduced for $\tilde{s}$ as long as 
$$z(1,0)=z(0,1)=1.$$

In particular, this implies that 
$$s(\lambda, 1)+s(0, 1)=0.$$

Furthermore, the condition
$$\partial_1 s(1,0)=1$$ 
is reproduced as long as 
$$\gmu=\partial_1 z (1,0).$$

We will now summarize the above discussion in the following definition of the renormalization operator acting on generating functions originally due  to the authors of \cite{EKW2}:
\begin{defn}\label{EKW_def}
\begin{eqnarray}
\label{ren_eq} {\cR}_{EKW}[s](x,y)&=&\gmu^{-1} s(z(x,y),\lambda y), {\rm where}\cr
\label{midpoint_eq} 0&=&s(\lambda x, z(x,y))+s(\lambda y, z(x,y)), \cr
\label{lambda-equation} 0&=&s(\lambda,1)+s(0,1) \quad {\rm and} \quad \gmu=\partial_1 z (1,0).
\end{eqnarray}
\end{defn}

\begin{defn}\label{B_space}
The Banach space of functions  $s(x,y)=\sum_{i,j=0}^{\infty}c_{i j} (x-\tau)^i (y-\tau)^j$, analytic on a bi-disk
$$ |x-\tau|<\rho, |y-\tau|<\rho,$$
for which the norm
$$\|s\|_\rho=\sum_{i,j=0}^{\infty}|c_{i j}|\rho^{i+j}$$
is finite, will be referred to as $\cA^{\tau}(\rho)$.

$\cA_s^{\tau}(\rho)$ will denote its symmetric subspace $\{s\in\cA^{\tau}(\rho) : s_1(x,y)=s_1(y,x)\}$. We will also use the shorthand notation 
$$\cA(\rho) \equiv \cA^0(\rho), \quad \cA_s(\rho) \equiv \cA^0_s(\rho).$$
\end{defn}

As we have already mentioned, the following  theorem has been proved with the help of a computer in \cite{EKW1} and \cite{EKW2}, and later in \cite{GJM}. We will quote here the version of the theorem from \cite{GJM}, since we will use the same functional spaces as the proofs in that paper.

\begin{oldthm}\label{GJMTheorem}

There exists a polynomial $s_a: \field{C}^2 \mapsto \field{C}$, such that 
\bigskip 
\begin{itemize}
\item[i)] The operator $\cR_{EKW}$ is well-defined, analytic and compact in $\cB_{r}(s_0) \subset \cA_s(\rho)$, with
$$\rho=1.75, \quad r= 1.1 \times 10^{-10}.$$

\bigskip 

\item[ii)] There exists a function $s^* \in \cB_r(s_0) \subset \cA_s(\rho)$ such that
$$\cR_{EKW}[s^*]=s^*.$$ 

\bigskip 

\item[iii)] The linear operator  $D \cR_{EKW}[s^{*}]$  has two eigenvalues outside of the unit circle:
$$ 8.72021484375 \le \delta_1 \le 8.72216796875, \quad  \delta_2={1 \over \lambda_*},$$
where 
$$  -0.248875313689    \le \lambda_* \le  -0.248886108398438.$$
\bigskip

\item[iv)] The complement of these two eigenvalues in the spectrum is compactly contained in the unit disk. The largest eigenvalue in the unit disk is equal to $\lambda_*$, while
%\begin{eqnarray}
%\nonumber &&{\rm spec}(D {\cR}_{EKW}[s^*]) \setminus \{\delta_1,\delta_2 \} \subset \{z \in \C: |z| \le |\lambda_*|\},\\
$$
{\rm spec}(D \cR_{EKW}[s^*]) \setminus \{\delta_1,\delta_2, \lambda_*\} \subset \{z \in \C: |z| \le 0.1258544921875\}.
$$
%\end{eqnarray}

\end{itemize}
 \end{oldthm}

%As we have already mentioned the following has been proved with the help of a computer in \cite{EKW}:
%\begin{thm}
%There is an $s^*$ in some Banach space of analytic functions, such that the operator ${\cR}_{EKW}$ is well-defined, analytic and compact on some neighborhood of $s^*$, and ${\cR}_{EKW}[s^*]=s^*.$ Furthermore, the scalings $\lambda^*$ and $\gmu^*$ corresponding to the fixed point $s^*$ satisfy
%\begin{eqnarray}
%\label{lambda} -0.2492 < &\lambda& < -0.2485, \\
%\label{mu} 0.0606 < &\gmu& < 0.0616.
%\end{eqnarray}
% \end{thm}

In this paper we will adopt a normalization condition different from $(\ref{lambda-equation})$ however. Specifically, $\mu$ will be defined from the normalization condition
$$
\cR_{EKW}[s](0,0)=1,
$$
i.e
\begin{eqnarray}\label{mu-equation}
\mu[s]= s(z(0,0),0). 
\end{eqnarray}

\begin{defn}
Define $\cI^\tau(\rho)$ and $\cI^\tau_s(\rho)$ to be the subsets of $\cA^\tau(\rho)$ and $\cA^\tau(\rho)$, respectively, of  functions normalized in the following way
\begin{eqnarray}
\nonumber \cI^\tau(\rho)&=&\left\{s \in  \cA^\tau(\rho): s(0,0)=1   \right\},\\
\nonumber \cI^\tau_s(\rho)&=&\left\{s \in  \cA^\tau_s(\rho): s(0,0)=1   \right\},\\
\nonumber \cI(\rho)&=&\left\{s \in  \cA(\rho): s(0,0)=1   \right\},\\
\nonumber \cI_s(\rho)&=&\left\{s \in  \cA_s(\rho): s(0,0)=1   \right\}.
\end{eqnarray}
\end{defn}

Our main result concerning the operator $\ref{EKW_def}$ will be the following

\begin{mainthmm}\label{analyticity_compactness}
There exists a polynomial $s_0 \in \cI(\rho)$, $\rho=1.75$, of  degree $9$
%\begin{eqnarray}
%\nonumber s_0(x,y)= 1\!&\!-\!&\!2.396026611328125 x -2.1419677734375 y \\
%\nonumber \!&\!+\!&\! 0.0618896484375  x^2   -2.13957113936142814 y^2 + 0.123779296875 x y \\
%\nonumber\! &\!-\!&\! 0.01153564453125 x^2 y -0.01092529296875x y^2,
%\end{eqnarray}
such that
\begin{itemize}
\item[i)] the operator $(\ref{EKW_def})$ is analytic in $B_{\delta}(s_0) \subset \cA(\rho)$ with $\delta=0.00405550003051758$;
\item[ii)] for all $s \in B_{\delta}(s_0)$ with real Taylor coefficients, the scalings  $\lambda=\lambda[s]$ and $\mu=\mu[s]$ satisfy
\begin{eqnarray}
\nonumber  0.000406771898269653 \le & \mu & \le 0.120654106140137,\\
\nonumber -0.276069164276123  \le  & \lambda & \le -0.222213745117188;
\end{eqnarray}  
\item[iii)] the operator $(\ref{EKW_def})$ is compact, with $\cR_{EKW}[s] \in \cI(\rho')$, $\rho'=1.0699462890625 \rho$, for all $s \in B_{\delta}(s_0)$.
\end{itemize}
\end{mainthmm}

\begin{rmk}
%Notice, that the choices of $\alpha$ and $\rho$ in the above Theorem are such that the fixed point of \cite{EKW2}:
%$$s^*(x,y)=s^*(x+\al,y+\al)=s^*((x+\al-{1 \over 2}) +{1\over 2},(y+\al-{1 \over 2}) +{1\over 2})$$ 
%is also in $\cA(\rho)$, $\rho=0.984375$. Indeed, this function is analytic on the set
%$$|x+\al-{1 \over 2})|<1.6, \quad |y+\al-{1 \over 2})|<1.6,$$
%which contains the bi-disk of radius $0.984375$ around $(0,0)$.

 We would like to note that the renormalization fixed point $s^*$ satisfies $$\|s^*-s_0\|_\rho  \le .00368565320968628,$$
(which is a rigorous bound), i.e., as expected, according to Theorem $\ref{analyticity_compactness}$, lies in the analyticity domain of renormalization.
\end{rmk}

\begin{rmk}
All the numbers quoted in the  Theorem and the Remark above are representable on a computer.
\end{rmk}

Although analyticity and compactness of the renormalization operator have been already proved in \cite{EKW2} (together with a stronger result of existence of the fixed point), we  would like to  reiterate that the goal of the paper will be to obtain these results with a lighter machinery than that of \cite{EKW2}. We also obtain a better bound on the size of the neighborhood of analyticity of $(\ref{EKW_def})$.

\section{Interval arithmetics in $\field{R}$ and $\field{R}^2$}

We will now give a very brief summary of interval operations in $\field{R}$ and $\field{R}^2$. For a more complete treatise of ``standard sets'' and operations on them, an interested reader is referred to an excellent review \cite{KSW}. 

A computer implementation of an arithmetic operation  $r_1 \# r_2$ ($\#$ is $+$, $-$, $*$ or $/$) on two real numbers does not generally yield an exact result. The ``computer'' result is a number representable in a standard IEEE floating point format (cf \cite{IEEE}). Such numbers are commonly referred to as ``representable''. In the $80$-bit extended precision IEEE arithmetics, a number is represented with $80$ bits of memory: $1$ bit for the sign of the number, $15$ bits for an exponents, and $64$ bits for a mantissa. A real representable number $r$ is of a the form 
$$r= \pm 1.m \ 2^{e-(2^{14}-1)},$$ 
where $\pm$ is chosen according to whether the sign bit is $0$ or $1$, $m$ is a base $2$ mantissa given  by the sequence of $0$'s and $1$'s associated with the $64$ bits of the mantissa, and the exponent $0 \le e \le 32767$ is defined by the state of the $15$ exponent bits (the state of $15$ ones is reserved to represent  ``overflows''). The representable  number $0$ is given by the  sequence of $80$ zeros. The set of all representable numbers will be denoted $\cR$. %Of specific interest are the smallest representable number larger than $1$, 
%\begin{equation}\label{u}
%u=1+2^{-99},
%\end{equation}
%and the largest representable number less than $1$, 
%\begin{equation}\label{d}
%d=1-2^{-99}.
%\end{equation}

Now, let $r_1 \# r_2$ be a mathematically legal, nonzero, arithmetic operation. The result of the true arithmetic operation  $r_1 \# r_2$ might not be a representable number. However, we can instruct the computer to attempt to round this operation either to the nearest representable number, up or down (our choice was to always round up). This might not be possible: the result of rounding up might have the exponent $e \ge  32767$ (overflow), or $e <0$ (underflow). In both cases the computer is instructed to raise an exception,  which is appropriately handled (either by terminating the program, or restarting with a different set of parameters). If, however, the result of rounding up is representable, the output of the computer implementation of the arithmetic operation is an upper bound on the true result of the operation. We will refer to such bound as ``the upper bound''.

We will define the {\it standard sets} in $\field{R}$ to be the collection of all closed real intervals $I[x,y]=\{r \in \field{R} | x \le r \le y  \}$:  
\begin{equation}
{\rm std}{(\field{R})}=\{I[x,y] \in \field{R}: x,y \in \cR \}.
\end{equation}  

If $I[x,y]$, $I[x_1,y_1]$ and $I[x_2,y_2]$ are in ${\rm std}{(\field{R})}$  then we can use the rounding up described above to obtain bounds on  the arithmetic operation on these sets as follows:

\begin{itemize}
\item[1)] {\it unary minus}:    $-(I[x,y])=I[-y,-x];$
\item[2)] {\it absolute value}: $|I[x,y]|=I[l,r]$, where  $l=\max\{ 0, x, -y \}$, $r=-\min\{0,x, -y\}$;
\item[3)] {\it addition}: $I[x_1,y_1]+ I[x_2,y_2]=I[x_3,y_3]$, where $y_3$ is the upper bound on $y_1+y_2$, while $-x_3$ is the upper bound on $-x_1+(-x_2)$;
\item[4)] {\it subtraction}: $I[x_1,y_1]- I[x_2,y_2] \equiv I[x_1,y_1]+ (-I[x_2,y_2])$;
\item[5)] {\it multiplication}: $I[x_1,y_1] \cdot I[x_2,y_2]=I[x_3,y_3]$, where $y_3$ is the maximum of the upper bounds on  $x_1  \cdot x_2$, $x_1 \cdot y_2$, $y_1  \cdot x_2$ and $y_1 \cdot y_2$, while $-x_3$ is the maximum of the upper bounds on $(-x_1)  \cdot x_2$, $(-x_1) \cdot y_2$, $(-y_1)  \cdot x_2$ and $(-y_1) \cdot y_2$; 
\item[6)] {\it inverse}: if $x_1 \cdot y_1 >  0$, then $I[x_1,y_1] \ {\rm inverse}=I[x_2,y_2]$, where $y_2$  is  the upper bound on $1/x_1$ and $-x_2$  is the  upper bound on $1/(-y_1)$;
%, ifd $x_1>0$ and $y_1>0$, and $y_2$  is  the upper bound on $1/x_1$ and $-x_2$  is the  upper bound on $1/(-y_1)$, ifd $x_1>0$ and $y_1>0$,
\item[7)] {\it division}:  if $x_2 \cdot y_2 >  0$, then $I[x_1,y_1] / I[x_2,y_2] \equiv I[x_1,y_1] \cdot (I[x_2,y_2] \ {\rm inverse})$. 
\end{itemize}

A standard set in $\field{R}^2$ is, naturally, a direct product of two standard sets in $\field{R}$:
\begin{equation}
{\rm std}{(\field{R}^2)}=\{I[x_1,y_1] \times I[x_2,y_2], x_1,x_2,y_1,y_2 \in \cR \}.
\end{equation}   

The arithmetic operation on these standard sets in $\field{R}^2$ are reducible to those on reals in an obvious way.

To obtain a bound on a algebraic and transcendental function, one can use their Taylor series together with a bound on the remainder. In our proofs we will require only three such functions: $\exp$, $\ln$ and $\sqrt{\phantom{a}}$. All of them are implemented via a finite truncation of their Taylor series with a bound on the remainder.

\section{Contraction Mapping Principle}

We will now outline a rather general method for finding a fixed point of a hyperbolic operator in a Banach space via its approximate Newton map. 

Let $C$ be an operator analytic and hyperbolic on some neighborhood $\cN$ in a Banach space $\cZ$. Suppose, that one knows its approximate hyperbolic fixed point  $Z_0 \in \cN$. Set
$$M \equiv \left[ \field{I}-D C[Z_0] \right]^{-1},$$
and for all $z$, such that $Z_0+M z \in \cN$,
$$
N[z]=z+C[Z_0+M z]-(Z_0+M z).$$

Notice, that if $z^*$ is a fixed point of $N$, then $Z_0+M z^*$ is a fixed point of $C$.

The linear operator $\field{I}-D C[Z_0]$ is indeed invertible since $D C$ is hyperbolic at $Z_0$. If $Z_0$ is a reasonably good approximation of the true fixed point of $C$, then the operator $N$ is expected to be a strong contraction in a neighborhood of $0$:
\begin{eqnarray}
\nonumber D N[z] &=& \field{I} +D C[Z_0+M z] \cdot M -M \\
\nonumber &=& \left[ M^{-1} +D C[Z_0+M z]  -\field{I} \right] \cdot M \\
\nonumber &=& \left[\field{I}-D C[Z_0]  +D C[Z_0+M z]  -\field{I} \right] \cdot M \\
\nonumber &=& \left[ D C[Z_0+M z] - D C[Z_0] \right] \cdot M.
\end{eqnarray}

The last expression is typically small in a small neighborhood of $0$, if the norm of $M$ is not too large (if $M$ is large, one might have to find a better approximation $Z_0$ and take a smaller neighborhood of $0$).  The following well-known Theorem specifies a sufficient condition for existence of the fixed point:

\begin{thm}({\it Contraction Mapping Principle})

Suppose that the operator $N$ is well-defined and analytic as a map from  $\cN \subset \cZ$ to $\cZ$, where $\cZ$ is some Banach space. Let $Z_0\in \cN$ and $B_\delta(Z_0) \subset \cN$ (an open ball of radius $\delta$ around $Z_0$) be such that
$$\| D N[Z]\| \le  \cD <1, $$
for any $Z \in  B_\delta(Z_0)$, and
$$\|N[Z_0]-Z_0\|\le \epsilon.$$

If $\epsilon < (1-\cD) \delta$ then the operator $N$ has a fixed point $Z^*$ in $B_\delta(Z_0)$, such that
$$\| Z^*-Z_0\| \le {\epsilon \over 1-\cD }.$$
\end{thm}

We will use the above Contraction Mapping Principle in several instances in our proofs below. In all those cases we will be verifying the hypothesis of the Contraction Mapping Principle using mild computer assistance.

\section{Analyticity of renormalization}
 
We consider the space $\cA(\rho)$, $\rho=1.75$. Define 
$$
s_0(x,y) \equiv D_0(x) y^3+A_0(x) y^2+B_0(x) y +C_0(x),
$$
a polynomial of degree $3$, with 
$$
D_0(x)= \sum_{i=0}^6 d_i x^i,\quad A_0(x)= \sum_{i=0}^6 a_i x^i, \quad B_0(x)= \sum_{i=0}^6 b_i x^i,\quad C_0(x)= \sum_{i=0}^6 c_i x^i,
$$
where the coefficients of these polynomials are as follows (the numbers given in the table are highly accurate approximations of the representable numbers actually used in the programs):

  \begin{center}
\begin{tabular}{|l|l|l|}
\hline
$i$ &  \hspace{2cm} $c_i$ &  \hspace{2cm} $b_i$  \\
\hline
$0$ & $\phantom{-} 1.00000000000000000$ & $-2.42962369607899157\times 10^{-1}$  \\
\hline
$1$ & $-1.02761956458970711$ & $\phantom{-} 5.87327440047455615\times 10^{-2}$ \\
\hline
$2$ & $\phantom{-} 2.93663720023727808\times 10^{-2}$ &
  $-5.93710236103475834\times 10^{-3}$  \\
\hline
$3$ &  $ -1.87658664952086400\times 10^{-3}$ & 
 $ \phantom{-} 6.09332694202817819\times 10^{-4}$ \\
\hline
$4$ &  $ \phantom{-}1.40668294317213841\times 10^{-4}$ &
  $ -6.46957663100331420\times 10^{-5}$  \\
\hline
$5$ &   $ -1.18664608613747513\times 10^{-5}$ & 
  $ \phantom{-} 7.02844653606969302\times 10^{-6}$  \\
\hline
$6$ &  $  \phantom{-}1.06935654680404746\times 10^{-6}$ &
  $  -7.75814237637266867\times 10^{-7}$ \\
\hline
\hline
$i$ &  \hspace{2cm} $a_i$ &  \hspace{2cm} $d_i$ \\
\hline
$0$ &  $-8.77647505670140721\times 10^{-1}$ & $-1.46791670728014469\times 10^{-3}$ \\
\hline
$1$ & $-5.62975994856259201\times 10^{-3}$ & $\phantom{-}5.62673177268855366\times 10^{-4}$ \\
\hline
$2$ & $\phantom{-}  9.13999041304226728\times 10^{-4}$&
 $-1.29391532620066284\times 10^{-4}$ \\
\hline
$3$ &
 $ -1.33251463600020727\times 10^{-4}$ &
 $ \phantom{-} 2.45371180652649971\times 10^{-5}$ \\
\hline
$4$ & 
  $\phantom{-} 1.84028385489487478\times 10^{-5}$ &
  $-4.18641657324405651\times 10^{-6}$ \\
\hline
$5$ & 
  $ -2.46417657757601024\times 10^{-6}$ &
  $\phantom{-}6.68319468409332288\times 10^{-7}$ \\
\hline
$6$ &
  $    \phantom{-}3.23459112039180563\times 10^{-7}$ &
  $  -1.01939399249366523\times 10^{-7}$. \\
\hline
\end{tabular}
\end{center}

%\begin{eqnarray}
%\nonumber a_0&=&-1.539804, \ a_1=-0.008287, \ a_2= 0.001262, \ a_3=-0.000172, \ a_4= 0.000022, \\
%\nonumber b_0&=&-1.276706, \ b_1= 0.092147, \ b_2=-0.008762, \ b_3= 0.000840, \ b_4=-0.000083, \\
%\nonumber c_0&=& 1.0,      \ c_1=-1.742031, \ c_2= 0.046273, \ c_3=-0.002762, \ c_4= 0.000193.
%\end{eqnarray}

A numerical study demonstrates that $s_0$ is a relatively good approximation of the renormalization fixed point  in $\cA_s(\rho)$ for a rather wide range of radii $\rho$ ($\rho \approx 0.8 \ldots 2.5$, below $\rho$ will be fixed to be $1.75$).

In this Section we will prove Theorem \ref{analyticity_compactness}. 
%We would like to note that the true renormalization fixed point $S^*$ satisfies $$\|S^*-S_0\|_r  \le 0.0065,$$ which is a rigorous bound, i.e., according to Theorem $\ref{analyticity_compactness}$, lies in the analyticity domain of renormalization.
This Theorem is proved with a mild aid of the computer. The computer assistance is restricted only to interval arithmetics on real numbers, (which, given sufficient time, in principle, can be done ``by hand''). 

We will start by noticing that a solution of the equation 
$$\label{Z-eqn}
s(x,Z')+s(y,Z')=0
$$
with $s \equiv s_0$ is given by
$$
\nonumber Z'_0(x,y) = {B_0(x)+B_0(y) \over 2 (A_0(x)+A_0(y) )} \left[\sqrt{1-4 {(A_0(x)+A_0(y) (\tilde{C}_0(x)+\tilde{C}_0(y) \over (B_0(x)+B_0(y))^2}}-1 \right],$$
where $\tilde{C}_0(x)$ is a bound on $C_0(x)+y^3 D(x)$, specifically:
\begin{equation}\label{Ctilde}
\tilde{C}_0(x)=C_0(x)+C  D_0(x), \quad C \in \field{D}_{\rho^3}(0).
\end{equation}
In the formula above we choose a branch of $\sqrt{}$ such that the range of $Z'_0$ is expected to be contained in $D_\rho$.

We would like to remark that we have intentionally chosen a polynomial $s_0$ of a sufficiently high degree so that it would be close to  the fixed point $s^*$  found in \cite{GJM}. This will later imply that the fixed point $s^*$ is indeed in the domain of analyticity and compactness of $\cR_{EKW}$. At the same time, taking a polynomial of degree higher than $3$ in $y$ is superfluous for the purposes of the computation of $Z'$. Therefore, these extra degrees (polynomial $D_0$) have been included as an ``error term''  in $(\ref{Ctilde})$. It is clear, that 
$$C_0(x)+y^3 D_0(x) \in \tilde{C}_0(x)$$
for all $(x,y)$ such that $|x|<\rho, |y|<\rho$.

We shall now check that for the given values of coefficients the function  $Z'_0$ is analytic on a bi-disk of radius $r=0.483119964599609$ centered at point $(0,0)$. %with 
%$$p \equiv {\al} (\lambda[s_0]-1).$$
Below $\cZ(r)$ will denote the Banach space of functions analytic on $\left\{(x,y):|x|<r, |y|<r\right\}$ equipped with the sup-norm, denoted by  $| \cdot |_r$. 
%Let 
%$$
%Z'_0(x,y)=Z''_0(x+p,y+p).
%$$
We have the following

\begin{lem}\label{Z0-range}
$Z'_0$ is in $\cZ(r)$ for $r=0.483119964599609$, and satisfies 
$$
|Z'_0|_r < 1.562789916992188.
$$ 
\end{lem}

\proof   First, define constants ${\ups}_0$, $\be_0$ and $\ga_0$, and quartic polynomials ${\ups}$, $\be$ and $\ga$, such that ${\ups}(0,0)=\be(0,0)=\ga(0,0)=0$ by setting
\begin{eqnarray}
\nonumber {\ups}_0+{\ups}(x,y)&=&A_0(x)+A_0(y), \\
\nonumber \be_0+\be(x,y)&=&B_0(x)+B_0(y), \\
\nonumber \ga_0+\ga(x,y)&=&\tilde{C}_0(x)+\tilde{C}_0(y).
\end{eqnarray}
Then
\begin{equation}\label{Z_0}
Z'_0(x,y)={\be_0+\be(x,y) \over 2 ({\ups}_0+{\ups}(x,y))} \left[\sqrt{1-4 {{\ups}_0 \ga_0 \over \be^2_0} + F_1(x,y)-F_2(x,y)}-1 \right],
\end{equation}
where
\begin{eqnarray}
\nonumber F_1(x,y)&=& {4{\ups}_0 \ga_0 \over \be^2_0 } { 2 \be_0 \be(x,y) +\be^2(x,y) \over \be^2_0+ 2 \be_0 \be(x,y) +\be^2(x,y)},\\
\nonumber F_2(X,Y)&=& 4 {{\ups}_0 \ga(x,y)+{\ups}(x,y) \ga_0+{\ups}(x,y) \ga(x,y) \over  \be^2_0+2 \be_0 \be(x,y) +\be^2(x,y) }.
\end{eqnarray}

Norms of both $F_1$ and $F_2$ are elementary:
$$ 
\nonumber |F_1|_r \le {4 |{\ups}_0 \ga_0| \over \be^2_0} {|2 \be_0 \be +\be^2|_r \over \be^2_0-|2 \be_0 \be +\be^2|_r} , \quad |F_2|_r \le  4 {|{\ups}_0 \ga+{\ups} \ga_0+{\ups} \ga|_r \over  \be^2_0-|2\be_0 \be + \be|^2_r}.
$$
Notice that all functions here whose norm has to be evaluated are explicit polynomials. We have estimated the above norms using the interval arithmetics on a computer. Next, we use these bounds in the following expression
$$
|Z'_0 |_r \le {|\be_0|+|\be|_r \over 2 (|{\ups}_0|-|{\ups}|_r) } \max_{j= \pm 1 }  \left\{  \sqrt{ 1- 4 {{\ups}_0 \ga_0 \over \be^2_0} + j |F_1|_r+j |F_2|_r }-1 \right\},
$$
and evaluate it on a computer.
\QED

Below, we will use the following shorthand notation for the supremum of $Z'_0$:
\begin{equation}\label{t}
t \equiv |Z'_0|_r.
\end{equation}

Notice that if one writes 
$$
s(x,y) \equiv \sigma(x,y) -{\tau(y) \over 2},
$$
where $\tau$ is chosen so that it has an inverse  branch $\nu$ with a range in $D_\rho$ (it is sufficient to choose $\tau$ a quadratic function), then a solution of equation $(\ref{Z-eqn})$ satisfies  
$$
Z'(x,y)=\nu[ \sigma(x,Z'(x,y))+\sigma(y,Z'(x,y)) ].
$$
(this ``trick'' has been shown to me by Hans Koch).

We will now specify the choices of $\tau$ and $\sigma$ for $s_0$. Set 
\begin{eqnarray}
\nonumber s_0(x,y) &\equiv& \sigma_0(x,y)-{\tau(y) \over 2}, \\
\nonumber \sigma_0(x,y) &\equiv& a(x) y^2+b(x)y +c(x), \quad {\rm and} \quad \Sigma_0(x,y)=\sigma_0(x+p,y),\\ 
\nonumber \tau(y) &\equiv& -2 a_0 y^2 -2 b_0 y, \\
\nonumber \nu (x) &\equiv& {-b_0 + \sqrt{b_0^2 -2 a_0 x} \over 2 a_0}.
\end{eqnarray}

It can be readily verified that $Z'_0(x,y)=\nu[\Sigma_0(x,Z'_0(x,y))+\Sigma_0(y,Z'_0(x,y))]$. We will next demonstrate that the operator $C_h$, defined on $\cZ(r)$ by setting
$$
 \nonumber C_h[Z'] \equiv \nu\left[  \Sigma_0 \circ ( \Pi_1 , Z' ) +h \circ (P_1,Z' )+\Sigma_0 \circ ( \Pi_2 , Z' ) +h \circ (P_2,Z' )  \right],
$$
where 
$$
P_1(x,y) \equiv x+p, \quad {\rm and} \quad P_2(x,y) \equiv y+p,
$$
has a fixed point for sufficiently small $h$'s.

\begin{prop}\label{Z-contract}
Let $\delta$, $\rho$ and $r$ be as in Theorem $\ref{analyticity_compactness}$ and Lemma $\ref{Z0-range}$. Then, for all  $s \in B_{\delta}(s_0) \subset \cA(\rho)$ the operator $C_{s-s_0}$ has a unique fixed point $Z'_s$ in $B_\epsilon(Z'_0) \subset  \cZ(r)$ with $\epsilon=0.01465$, and the map $s \mapsto Z'_s$ is analytic from  $B_{\delta}(s_0)$ to $B_\epsilon(Z'_0)$.
\end{prop}
\proof Define an operator 
$$
N_h[z]=z+C_h[Z'_0+M z]-(Z'_0+M z), \quad M \equiv \left[ \field{I}-D C_0[Z'_0] \right]^{-1}.
$$
We will demonstrate that this ``Newton map'' has a fixed point in a neighborhood of $0$.

We will first estimate the norm of $N_h[0]=C_h[Z'_0]-Z'_0$:
\begin{eqnarray}\label{C_diff}
\nonumber |C_h[Z'_0]-Z'_0|_r&=&\left| \nu\left[  \Sigma_0 \circ (\Pi_1, Z'_0) +h\circ (\Pi_1,Z'_0)+ \Sigma_0\circ (\Pi_2,Z'_0) +h\circ (\Pi_2,Z'_0) \right]\right. \\
\nonumber  & & \phantom{aa} -\left. \nu\left[  \Sigma_0(\Pi_1, Z'_0) + \Sigma_0(\Pi_2,Z'_0) \right] \right|_r\\
\nonumber &=& {1 \over 2 a_0} \left| \sqrt{b_0^2 -2 a_0 \left[ \Sigma_0(\Pi_1, Z'_0) +h(P_1,Z'_0)+ \Sigma_0(\Pi_2,Z'_0) +h(P_2,Z'_0) \right]} \right.\\
\nonumber & & \phantom{aaaa}- \left. \sqrt{b_0^2 -2 a_0 \left[ \Sigma_0(\Pi_1, Z'_0)+ \Sigma_0(\Pi_2,Z'_0) \right]} \right|_r \\
\nonumber &=& {1 \over 2 a_0} \left| F(1)-F(0) \right|_r,
\end{eqnarray}
where 
$$
F(\zeta)= \sqrt{b_0^2 -2 a_0 \left[ \Sigma_0(\Pi_1, Z'_0) +\zeta h(P_1,Z'_0)+ \Sigma_0(\Pi_2,Z'_0) +\zeta h(P_2,Z'_0) \right]}
$$
is an analytic function from $\field{D}_R$ to $\cZ(r)$ (here and below, $\field{D}_R$ stands for an open disk of radius $R$ in $\field{C}$), with 
\begin{equation}\label{R}
R \le { b_0^2-4 a_0 \theta_0  -2 |a_0| | \Theta_0(\Pi_1, Z'_0)|_r -2 |a_0| |\Theta_0(\Pi_2,Z'_0)|_r \over 2 |a_0| |h(P_1,Z'_0)|_r + 2 |a_0| |h(P_2,Z'_0)|_r },
\end{equation}
where 
$$
\nonumber \theta_0 \equiv \Sigma_0(0,0),\quad \Theta_0(x,y) \equiv \Sigma_0(x,y)-\theta_0.
$$
Below, we will denote the (positive) real number $b_0^2-4 a_0 \theta_0$ by $c$.  
A straightforward Cauchy estimate yields:
$$
\nonumber \left| F(1)-F(0) \right| \le {1 \over R-1} \sup_{|\zeta|\le R} |F(\zeta)| \le { 4  |a_0| \delta  \sqrt{2 c} \over  c  -4  |a_0| | \Theta_0|_{(r,t)}  - 4  |a_0| \delta}.
$$
%where $\varrho=\max\{r+|p|, t \}$, and $t$ is as in $(\ref{t})$. 
We will denote 
\begin{equation}\label{M_diff}
\varepsilon \equiv{ 2  \delta  \sqrt{2 c} \over   c -4  |a_0| | \Theta_0|_{(r,t)} - 4 |a_0| \delta}.
\end{equation}

This expression has been shown to satisfy
$$
\varepsilon < 0.0137615203857422.
$$

As a next step we will provide a bound on the derivative of the operator $N_h$:
\begin{eqnarray}
\nonumber | D N_h[z] | &=& \left|\left[ D C_h[Z'_0+M z]-D C_0[Z'_0] \right] M\right| \\
\nonumber &\le& \left| \left[ D C_h[Z'_0+M z]-D C_0[Z'_0] \right] \right| |M|.
\end{eqnarray}

To estimate the norm of $M$ we first verify that  $|D C_0[Z'_0]|<1$:
\begin{eqnarray}
\nonumber \left| D C_0[Z'_0]\right|&=&\left| {\partial_2 \Theta_0(\Pi_1, Z'_0)+ \partial_2 \Theta_0(\Pi_2,Z'_0) \over 2 \sqrt{c -2 a \left[ \Theta_0(\Pi_1, Z'_0)+ \Theta_0(\Pi_2,Z'_0) \right] }} \right| \\
\nonumber &\le& {\left| \partial_2 \Theta_0 \right|_{(r,t)}  \over  \sqrt{c - 4 |a_0| \left| \Theta_0 \right|_{(r,t)} }}.
\end{eqnarray}
A computer-aided evaluation of the last expression indeed shows that it is less than $1$. This implies that $M$ can be found as the limit of a convergent series whose norm is easily bounded:
\begin{equation}\label{M}
| M| \le {1 \over 1-\left| D C_0[Z'_0]\right|}< 1.0430755615234375 \equiv \cM.
\end{equation}

Next, 
\begin{eqnarray} \label{D_diff}
 \nonumber  D C_h[Z']-D C_0[Z'_0]& = &{\partial_2 \Theta_0(\Pi_1, Z')+ \partial_2 \Theta_0(\Pi_2,Z')+\partial_2 h(P_1, Z')+ \partial_2 h(P_2,Z') \over 2 \sqrt{c -2 a_0 \left[ \Theta_0(\Pi_1, Z')+\Theta_0(\Pi_2,Z')+  h(P_1,Z')+  h(P_2,Z')  \right] }}\\
\nonumber  &-&{\partial_2 \Theta_0(\Pi_1, Z'_0)+ \partial_2 \Theta_0(\Pi_2,Z'_0) \over 2 \sqrt{c -2 a_0 \left[ \Theta_0(\Pi_1, Z'_0)+ \Theta_0(\Pi_2,Z'_0) \right] }}\\
 \nonumber  & = &\left[ {\partial_2 \Theta_0(\Pi_1, Z')+ \partial_2 \Theta_0(\Pi_2,Z')+\partial_2 h(P_1, Z')+ \partial_2 h(P_2,Z') \over 2 \sqrt{c -2 a_0 \left[ \Theta_0(\Pi_1, Z')+  \Theta_0(\Pi_2,Z')+  h(P_1,Z')+  h(P_2,Z')  \right] }} \right.\\
\nonumber  &-&\left. {\partial_2 \Theta_0(\Pi_1, Z')+ \partial_2 \Theta_0(\Pi_2,Z') \over 2 \sqrt{c -2 a_0 \left[ \Theta_0(\Pi_1, Z')+ \Theta_0(\Pi_2,Z')+  h(P_1,Z')+ h(P_2,Z')   \right] }} \right]\\
\nonumber &+& \left[ {\partial_2 \Theta_0(\Pi_1, Z')+ \partial_2 \Theta_0(\Pi_2,Z') \over 2 \sqrt{c -2 a_0 \left[ \Theta_0(\Pi_1, Z')+ \Theta_0(\Pi_2,Z')+  h(P_1,Z')+  h(P_2,Z')   \right] }} \right.\\
\nonumber  &-& \left.{\partial_2 \Theta_0(\Pi_1, Z')+ \partial_2 \Theta_0(\Pi_2,Z') \over 2 \sqrt{c -2 a_0 \left[ \Theta_0(\Pi_1, Z')+ \Theta_0(\Pi_2,Z')  \right] }}\right] \\
\nonumber &+& \left[{\partial_2 \Theta_0(\Pi_1, Z')+ \partial_2 \Theta_0(\Pi_2,Z') \over 2 \sqrt{c -2 a_0 \left[ \Theta_0(\Pi_1, Z')+ \Theta_0(\Pi_2,Z')  \right] }}\right. \\
\nonumber &-&\left.{\partial_2 \Theta_0(\Pi_1, Z'_0)+ \partial_2 \Theta_0(\Pi_2,Z'_0) \over 2 \sqrt{c -2 a_0 \left[ \Theta_0(\Pi_1, Z'_0)+  \Theta_0(\Pi_2,Z'_0) \right] }}\right]\\
 &\equiv& I_1+I_2+I_3.
\end{eqnarray}
We will estimate norms of the three expressions in brackets separately.
\begin{equation}\label{I1}
|I_1|_r \le {  |\partial_2 h|_{(r,s)} \over \sqrt{c -4 |a_0| |\Theta_0|_{(r,s)} -4 |a_0|  |h|_{(r,s)} } } \le  { m \delta  \over \sqrt{ c -4 |a_0| |\Theta_0|_{(r,s)} -4 |a_0|   \delta  }  },
\end{equation}
where
%$$
%s=t + |M| \epsilon, \quad  m=n {s^{n-1} \over \rho^n }, \quad n={1 \over \ln{\left[ { \rho \over s} \right]}}, \quad {\rm and } \quad \varrho=\max\{r+|p|, s \}.
%$$
$$
s=t + |M| \epsilon, \quad  m=n {s^{n-1} \over \rho^n } \quad {\rm and} \quad  n={1 \over \ln{\left[ { \rho \over s} \right]}}.
$$
We use a Cauchy estimate for the function 
%$$
$$F_2(\zeta)={1 \over \sqrt{ c -2 a_0 \left[ \Theta_0(\Pi_1, Z')+ \Theta_0(\Pi_2,Z')+  \zeta h(P_1,Z')+  \zeta h(P_2,Z')   \right] }}$$
%$$
to bound $|I_2|_r$. Notice, that, $F_2$ is also analytic on $\field{D}_R$, but gets unbounded as $\zeta$ approaches the boundary of this disk. Therefore, to bound $|I_2|_r$, we use   a Cauchy bound  for $F_2$ on a smaller disk of radius
%$$\tilde{R}={{2 \over 3} c  -2 |a_0| | \Theta_0(\Pi_1, Z'_0)|_r -2 |a_0| |\Theta_0(\Pi_2,Z'_0)|_r \over 2 |a_0| |h(P_1,Z'_0)|_r + 2 |a_0| |h(P_2,Z'_0)|_r }:$$

$$\tilde{R}={c  -4 |a_0|  \left| \Theta_0 \right|_{(r,s)} \over 4 |a_0| \delta } - {1 \over 3} {c  -4 |a_0|  \left| \Theta_0 \right|_{(r,s)} - 4 |a_0| \delta  \over 4 |a_0| \delta}:$$

$$
|I_2|_r \le \left| \partial_2 \Theta_0 \right|_{(r,s)} {1 \over \tilde{R}-1}  \sup_{|\zeta|\le \tilde{R}} |F_2(\zeta)| \le  { 4 |a_0|  \delta    \left| \partial_2 \Theta_0 \right|_{(r,s)}  \over 2 \left({1 \over 3} \left({c  -4 |a_0|  \left| \Theta_0 \right|_{(r,s)} - 4 |a_0| \delta } \right)\right)^{3/2} }.
$$
Finally, to estimate $|I_3|_r$, we use a Cauchy bound for the function 
$$
F_3(\zeta)={\partial_2 \Theta_0(\Pi_1, Z'_0 +\zeta M z)+ \partial_2 \Theta_0(\Pi_2,Z'_0+\zeta M z)  \over 2 \sqrt{ c -2 a_0 \left[ \Theta_0(\Pi_1, Z'_0+\zeta M z)+ \Theta_0(\Pi_2,Z'_0 + \zeta M z)   \right] }},
$$
analytic on $\field{D}_{\hat{R}}$ with 
$$
\hat{R}={\rho -t \over |M| \epsilon}:
$$
\begin{equation}\label{I3}
|I_3|_r=|F_3(1)-F_3(0)| \le {1 \over \hat{R}-1} \sup_{|\zeta| \le \hat{R}} |F_3(\zeta)| \le {|M| \epsilon \over \rho-s}  {|\partial_2  \Theta_0 |_{(r,\rho)} \over \sqrt{ c -4 |a_0| | \Theta_0|_{(r,\rho)} }}.
 \end{equation}

Individual norms $|M|$, $|I_1|_r$, $|I_2|_r$ and $|I_3|_r$ and their sum have been estimated on a computer to produce:
\begin{equation} \label{M_der}
|D N_h(z)| \le  0.0125999450683594 \equiv \cD
\end{equation}
for all $\|h\|_\rho < \delta$ and $|z|_r < \epsilon$.

Estimates $(\ref{M_diff})$ and $(\ref{M_der})$ demonstrate that 
$$\varepsilon < (1-\cD) {\epsilon \over \cM},$$
 and the claim follows from the Contraction Mapping Principle.
\QED

At the next step we obtain bounds on the scaling parameters $\lambda$ and $\mu$.

Recall that the scaling $\lambda$ solves equation $(\ref{lambda-equation})$.

The prove Theorem $\ref{analyticity_compactness}$ we will require the following

\begin{lem}\label{lambda_lemma}
Equation $(\ref{lambda-equation})$ has a solution $\lambda=\lambda[s]$ for all $s \in   B_{\delta}(s_0)$, and satisfies 
$$
-0.276069164276123  \le   \lambda  \le -0.222213745117188
$$
for all $s \in B_\delta (s_0)$ with real Taylor coefficients.

Furthermore, the map $s \mapsto \lambda[s]$ is analytic on $B_\delta (s_0)$.
\end{lem}
\proof 
Write,
$$s_1(x,y) \equiv A_1(y)x^2 +B_1(y)x +C_1(y),
$$
where 
\begin{eqnarray}
\nonumber A_1(y)&=&c_2+b_2 y+a_2 y^2+ d_2 y^3, \\
\nonumber B_1(y)&=&c_1+b_1 y+a_1 y^2+ d_1 y^3, \\
\nonumber C_1(y)&=&c_0+b_0 y+a_0 y^2+ d_0 y^3.
\end{eqnarray}

Denote

Notice, that  the equation $(\ref{lambda-equation})$ with $s=s_1$ can be solved exactly (a quadratic equation for $\lambda$). Furthermore, suppose $s=s_1+h$, $\| h\|_\rho \le  \|s_0-s_1\|_\rho+ \delta \equiv \tilde{\delta}$, then $\lambda[s]$ solves 
\begin{equation}\label{lambdaS-equation}
s_1(\lambda[s],1)+s_1(0,1)=- h(\lambda[s],1)-h(0,1).
\end{equation}

Write 
$$\lambda[s]=\lambda[s_1]+\Delta \lambda,
$$
then $\Delta \lambda$ is a solution of 
\begin{equation}\label{lambdaD-equation}
A_1(1) \cdot (\Delta \lambda)^2+[2 A_1(1) \cdot  \lambda[s_1]+B_1(1)] \cdot \Delta \lambda+ C=0
\end{equation}
with $C \in \field{D}_{2 \tilde{\delta}}$ (an open disk of radius $2 \tilde{\delta}$ in $\field{C}$).
% where
%\begin{eqnarray}
% \nonumber \varrho&=&{\rho_1 +\rho_2 \over \rho}  \delta,\\      
% \nonumber \rho_1&=&\max\{|{\al}|,|1-{\al}|\},\\
%\nonumber  \rho_2&=&\max\{R,|1-{\al}| \},
%\end{eqnarray}
%and $R$ is some {\it a priori} bound on $|\lambda[s]-{\al}|$. 
For all $s$ with real coefficients, (\ref{lambdaD-equation}) is a quadratic equation for $\Delta \lambda$ with real coefficients, whose solution gives the required bound.

%We verify {\it a posteriori} that 
%$$
%R>|\lambda[s_0]-{\al}|+\sup{|\Delta \lambda|}.
%$$
\QED

The last ingredient in the proof of analyticity of renormalization is a bound on the scaling $\mu$.

\begin{lem}\label{mu_lemma}
Let $\mu[s]$ be as in equation $(\ref{mu-equation})$. Then the map $s \mapsto \mu[s]$ is analytic on $B_\delta (s_0)$, and $\mu[s]$ satisfies
$$
0.000406771898269653 \le  \mu  \le 0.120654106140137,
$$
for $s \in B_\delta (s_0)$ with real Taylor coefficients.
\end{lem}
\proof  The bound is straightforward:
\begin{eqnarray}
\nonumber \mu[s]&=&s(Z'_s(0,0),0)\\
\nonumber           & \in &s_0(Z'_s(0,0),0)+h(Z'_s(0,0), 0). 
\end{eqnarray}

Now, it is  clear, that whenever $s$ has real coefficients,
$$
\mu[s] \in s_0(Z'_0(0,0)+I[-|M| \epsilon, |M| \epsilon],0)+I\left[-\delta,\delta \right], 
$$
where, as before, $I[x,y]$ stands for a closed real interval with endpoints $x$ and $y$,
%, and
%$$
%R=\max\{t+|M| \epsilon, {\al}|\lambda[s]-1| \}, 
%$$
while the rest of the quantities are as in Lemma $\ref{lambda_lemma}$ and Prop. $\ref{Z-contract}$.
\QED

Now, part $i)$ of  Theorem $\ref{analyticity_compactness}$ follows immediately from  Proposition \ref{Z-contract} and Lemmas \ref{lambda_lemma} and \ref{mu_lemma}.

\section{Compactness of renormalization}

We will now outline the proof of compactness of the renormalization operator.

\begin{prop}\label{compactess}
For every $s \in B_\delta(s_0) \subset \cA(\rho)$ its renormalization $\cR_{EKW}[s]$ is in $\cI(\kappa \rho)$ with $\kappa=1.0699996948242188.$.
\end{prop}

\proof Let  $r'=\kappa r$ and $\rho'=\kappa \rho$. Verification that $|\lambda[s]| \rho'| \le \rho$ is straightforward. The ``difficult'' part is to demonstrate that $|Z'_s|_{r'} \le \rho$. To show this, we will use the fact that $Z'_s=\lim_{n\rightarrow \infty}N_{s-s_0}[0]$,   

To that end we first verify that $Z'_0 \in \cZ(r')$. This amounts to checking on a computer that 
$$
\nonumber |{\ups}_0| > |{\ups}|_{r'} \quad {\rm and} \quad 1-4 {{\ups}_0 \ga_0 \over 
\be^2_0} > | F_1|_{r'}+|F_2|_{r'} \quad {\rm (cf.} \quad (\ref{Z_0})).
$$

Below we will use the following shorthand notation:  $t' \equiv |Z'_0|_{r'}$.

Set $z_0=0$, and define  for all $n \ge 0$
$$
z_{n+1} \equiv  N_h[z_{n}], \quad \Delta z_{n}=z_{n+1}-z_n, \quad  Z'_n=Z'_0+M z_n
$$
 (we have suppressed the dependence on $s$ in $Z'_n$ for notational convenience), then  
\begin{eqnarray}
\nonumber |\Delta z_0|_{r'}&=& |C_h[Z'_0]-Z'_0|_{r'} \\
\nonumber &\le& { 2  \delta  \sqrt{2 c} \over  c -4 |a_0| | \Theta_0|_{(r',t')} - 4 |a_0| \delta}  \quad {\rm (cf. \quad equation \quad (\ref{M_diff}) )},
\end{eqnarray}
%where 
%$\varrho \equiv \max \{r'+|p|,t'  \}$ (we also verify that $\varrho<\rho$), 
and
\begin{eqnarray}
\nonumber |\Delta z_n|_{r'}&=&|N_h[z_n]-N_h[z_{n-1}]|_{r'}\\
\nonumber &=&| \Delta  z_{n-1}+C_h[Z'_{n-1}+M \Delta z_{n-1}]-C_h[Z'_{n-1}] -M \Delta 
z_{n-1}|_{r'}\\
\nonumber & \le &| \left[C_h[Z'_{n-1}+M \Delta z_{n-1}]-C_h[Z'_{n-1}]-D C_h[Z'_{n-1}] M 
\Delta z_{n-1}\right]|_{r'}+\\
\nonumber &\phantom{=}&|\left[D C_h[Z'_{n-1}]-D C_h[Z'_0]  \right] M \Delta 
z_{n-1}|_{r'}\\
\nonumber &\equiv& T_1+T_2.
\end{eqnarray}
To evaluate $T_2$ we will require estimates similar to $(\ref{M})$--$(\ref{I3})$. In fact, 
estimates $(\ref{M})$--$(\ref{I3})$ can be 
recycled after one substitutes $r'$ for $r$, $\epsilon'$ for $\epsilon$  and $s' \equiv t'+|M| \epsilon'$ for $s$ (here $\epsilon'$ is an {\it a priori} bound on $|z_n|_{r'}$, $n \ge0$,  again,  verifiable {\it a posteriori}). We will not repeat the details, but rather use the following symbolic notation  
$$
|T_2|_{r'} \le C_2 |\Delta z_{n-1}|_{r'},
$$
keeping in mind that $C_2$ is computable on a machine.

To estimate $T_1$, we use a Cauchy bound for the following function:
\begin{eqnarray}
\nonumber F(\zeta)&=& \sqrt{c \!-\!2 a_0 \left[ \Theta_0(\Pi_1, Z'_{n-1}\!+\!\zeta M 
\Delta z_{n-1}) \!+\!h(P_1,Z'_{n-1}\!+\!\zeta M \Delta z_{n-1})\right.+}\\
\nonumber & &\\
\nonumber & &\overline{\left.\phantom{aaaaaaaa}\Theta_0(\Pi_2,Z'_{n-1}\!+\!\zeta M 
\Delta 
z_{n-1}) \!+\!h(P_2,Z'_{n-1}\!+\!\zeta M \Delta z_{n-1}) \right]},
\end{eqnarray}
analytic on $\field{D}_R$ with 
$$
R \le { \rho -|Z'_{n-1}|_{r'} \over |M| |\Delta z_{n-1}|_{r'} }.
$$
This results in
\begin{eqnarray}
\nonumber |T_1|_{r'} &\le& {1 \over R(R-1)}{1 \over 2 |a_0|} \sup_{|\zeta|\le R} 
|F(\zeta)|_{r'}\\
\nonumber &=&{|M|^2 |\Delta z_{n-1}|^2_{r'} \over 2 |a_0| (\rho-s') (\rho-s'-|M||\Delta z_{n-1}|_{r'})} \sqrt{c +4 |a_0| \left[  
|\Theta_0|_{(r',s')} +\delta \right]} \\
&\equiv& C_1 {|\Delta z_{n-1}|^2_{r'} \over \rho-s'-|M||\Delta z_{n-1}|_{r'}}.
\end{eqnarray}
%where $\varrho=\max\{r'+|p|, s' \}$. 
Finally, we have 
\begin{eqnarray}
\nonumber |\Delta z_{n}|_{r'} &\le&  C_1 {|\Delta z_{n-1}|^2_{r'} \over \rho-s'-|M||\Delta z_{n-1}|_{r'}} + C_2 |\Delta z_{n-1}|_{r'} \\
\nonumber &=& K_{n-1} |\Delta z_{n-1}|_{r'}|,\\
 \nonumber  K_{n-1} & \equiv & C_1 {|\Delta z_{n-1}|_{r'} \over \rho-s'-|M||\Delta z_{n-1}|_{r'}} + C_2
\end{eqnarray}
(we  check that $K_0<1$).

Clearly, $K_n$ is a monotone increasing function of $|\Delta z_{n}|_{r'}$, therefore $K_n<K_{n-1}$. We can now verify that the  sum
$$
\Sigma \equiv \sum^{\infty}_{k=0} |\Delta z_{k}|_{r'} 
$$

is convergent:
\begin{eqnarray}
\nonumber \Sigma &\le&  \sum^{\infty}_{k=1} K_{k-1}  |\Delta z_{k-1}|_{r'} +|\Delta z_0|_{r'} \\
\nonumber  &\le&  \sum^{\infty}_{k=2} K_{k-1} K_{k-2}  |\Delta z_{k-2}|_{r'} +K_0 |\Delta z_0|_{r'} +  |\Delta z_0|_{r'} \\
\nonumber &\le& \sum^{\infty}_{k=0} K_0^k  |\Delta z_{0}|_{r'} \\
\nonumber & \le&  {|\Delta z_{0}|_{r'} \over 1-K_0}.
\end{eqnarray}

Therefore, for all $n \ge 0$ 
$$
|Z'_n-Z'_0|_{r'} \le |M|  {|\Delta z_{0}|_{r'} \over 1-K_0} \equiv \gamma'.
$$
Finally, we verify that 
$$t'+\gamma'<\rho.$$

\QED

\section{Acknowledgments}

The author is thankful to Hans Koch for many useful discussions  on the subject.

\end{document}